\documentclass[final,1p,times]{elsarticle}
\usepackage{amsthm}
\usepackage{amsmath}
\usepackage{amssymb}
\usepackage{amsfonts}
\usepackage{graphicx}
\usepackage{color}
\usepackage{multirow}
\usepackage{pdflscape}

\usepackage{lipsum}
\makeatletter
\def\ps@pprintTitle{%
 \let\@oddhead\@empty
 \let\@evenhead\@empty
 \def\@oddfoot{}%
 \let\@evenfoot\@oddfoot}
\makeatother

\biboptions{sort&compress}

\begin{document}

\begin{frontmatter}

\title{Discrete Ricci curvatures for directed networks}

\author[label1]{Emil Saucan}
\author[label2]{R.P. Sreejith}
\author[label2]{R.P. Vivek-Ananth}
\author[label3,label4]{J\"urgen Jost}
\author[label2,label3]{Areejit Samal}
\address[label1]{Department of Applied Mathematics, ORT Braude College, Karmiel 2161002 Israel}
\address[label2]{The Institute of Mathematical Sciences (IMSc), Homi Bhabha National Institute (HBNI), Chennai 600113 India}
\address[label3]{Max Planck Institute for Mathematics in the Sciences, Leipzig 04103 Germany}
\address[label4]{The Santa Fe Institute, Santa Fe, New Mexico 87501 USA}
\cortext[cor1]{asamal@imsc.res.in}

\begin{abstract}
A goal in network science is the geometrical characterization of complex networks. In this direction, we have recently introduced Forman's discretization of Ricci curvature to the realm of undirected networks. Investigation of this edge-centric network measure, Forman-Ricci curvature, in diverse model and real-world undirected networks revealed that the curvature measure captures several aspects of the organization of undirected complex networks. However, many important real-world networks are inherently directed in nature, and the definition of the Forman-Ricci curvature for undirected networks is unsuitable for the analysis of such directed networks. Hence, we here extend the Forman-Ricci curvature for undirected networks to the case of directed networks. The simple mathematical formula for the Forman-Ricci curvature of a directed edge elegantly incorporates vertex weights, edge weights and edge direction. Furthermore we have compared the Forman-Ricci curvature with the adaptation to directed networks of another discrete notion of Ricci curvature, namely, the well established Ollivier-Ricci curvature. However, the two above-mentioned curvature measures do not account for higher-order correlations between vertices. To this end, we adjusted Forman's original definition of Ricci curvature to account for directed simplicial complexes and also explored the potential of this new, augmented type of Forman-Ricci curvature, in directed complex networks.
\end{abstract}

\end{frontmatter}

\section{Introduction}
\label{introduction}

In the last two decades, network theory \cite{Watts1998,Barabasi1999,Albert2002,Newman2010,Fortunato2010,Dorogovtsev2013} has become a standard tool to investigate complex networks that pervade all aspects and strata of nature and society. To wit, metabolic networks \cite{Jeong2000} capture interactions among different metabolites and enzymes that are responsible for growth and maintenance of a living organism. Ecological networks \cite{Sole2001,Montoya2006} model interactions among different species in the ecosystem. Online social networks such as Facebook \cite{Ellison2007} and Twitter \cite{Mcauley2012} encapsulate relationships between different individuals in the society. Transportation networks \cite{Subelj2011} encapture the movement of traffic across various parts of the world. Network science \cite{Watts1998,Barabasi1999,Albert2002,Newman2010,Fortunato2010,Dorogovtsev2013} aims to characterize the structure of these ubiquitous complex networks. Towards this goal, a recent focus in network science has been the development of geometry inspired measures to characterize the structure of complex networks \cite{Eckmann2002,Ollivier2009,Lin2010,Lin2011,Bauer2012,Jost2014,Wu2015,Ni2015,Sandhu2015,Sreejith2016,Bianconi2017,Samal2018,Courtney2018}.

Curvature represents a central concept in geometry, and it certainly is the focal point of differential and Riemannian geometry, which formally quantifies our intuition of the deviation of an object from being \textit{flat}. In the classical context, where the notion of curvature originated, several types of curvature have been considered \cite{Jost2017}. Among the different notions of curvature, it appears that the concept of Ricci curvature is the most useful, after a proper discretization process, for the analysis of graphs or networks \cite{Ollivier2009,Lin2010,Lin2011,Bauer2012,Jost2014,Ni2015,Sandhu2015,Sreejith2016,Samal2018}. To this end, let us note that, in differential geometry, Ricci curvature first appears in the context of the so called Jacobi equation, that represents a linearization of the equation for geodesics \cite{Jost2017} and, as such, governs both the growth of volumes and the dispersion rate of geodesics \cite{Lott2009,Jost2017}. Furthermore, Ricci curvature is also an essential ingredient in the so called Bochner-Weitzenb\"{o}ck formula, that connects between the Laplacian on a manifold and its geometry, as expressed by its curvature. By its very definition through the above mentioned equation of geodesics, Ricci curvature depends on a direction, which then for graphs or networks translates into the fact that it should be a measure associated to an edge rather than to a vertex. Among the several curvature measures \cite{Eckmann2002,Shavitt2004,Saucan2005,Ollivier2009,Ollivier2010,Lin2010,Lin2011,Narayan2011,Bauer2012,Jost2014,Ollivier2013,Wu2015,Ni2015,Sandhu2015,Sreejith2016} that have been proposed for geometrical characterization of complex networks, two different discretizations of the Ricci curvature, Ollivier-Ricci curvature \cite{Ollivier2009,Ollivier2010,Ollivier2013} and Forman-Ricci curvature \cite{Forman2003} appear particularly appealing and useful for analyses of complex networks. It is important to remark that Ollivier-Ricci curvature directly captures, in a discrete setting, the growth of volumes property of Ricci curvature. In contrast, Forman-Ricci curvature stems from a discretization of the Bochner-Weitzenb\"{o}ck formula which encodes, in a discrete setting, the dispersion of geodesics property of Ricci curvature.

Among the two discretizations of the Ricci curvature, Ollivier-Ricci curvature was introduced first to the realm of complex networks, and this notion has been systematically explored in both model and real-world undirected networks \cite{Lin2010,Lin2011,Bauer2012,Jost2014,Loisel2014,Ni2015,Sandhu2015,Sandhu2016}. More recently \cite{Sreejith2016}, we had adapted Forman's discretization of the classical Ricci curvature \cite{Forman2003} to the realm of undirected complex networks. Importantly, Forman-Ricci curvature of an edge due to its enticingly simple mathematical formula is eminently suitable for the analysis of large complex networks \cite{Sreejith2016}. Vitally, the mathematical formula of the Forman-Ricci curvature for an edge also elegantly incorporates both edge weights and vertex weights, and this also makes the measure amenable to the analysis of both unweighted and weighted networks \cite{Sreejith2016}. Since both Ollivier-Ricci and Forman-Ricci curvatures represent discretizations of the classical Ricci curvature which is intrinsically associated with edges of a network, the two notions of curvature does not necessitate the technical artifice of extending a measure for the curvature of vertices to the edges \cite{Sreejith2016,Samal2018}. Thus, both Ollivier-Ricci and Forman-Ricci curvatures can be exploited for edge-based analysis of complex networks \cite{Samal2018}. Moreover, one can easily extend the definitions of the Ollivier-Ricci and Forman-Ricci curvatures for an edge to define corresponding  Ollivier-Ricci and Forman-Ricci curvatures for a vertex \cite{Sandhu2015,Ni2015,Sreejith2016,Sreejith2017} by summing the curvatures of edges adjacent to a vertex in an undirected network, and this definition of the vertex curvature based on edge curvatures is in effect a discrete analogue of scalar curvature in the classical context \cite{Jost2017}. Specific to Forman-Ricci curvature, we remark that though the mathematical expression of the Forman-Ricci curvature for an edge is a local measure dependent on the weights of adjacent vertices and edges in the network \cite{Sreejith2016}, still, the local geometric characteristic can provide deep, nontrivial insights on the global topology of the network \cite{Forman2003} that go beyond the simple combinatorial structure of the network. Note that it is easy to construct simple examples of undirected networks which have the same degree distribution but very different distributions of the two notions of discrete Ricci curvature.

It is imperative to notice that several important real networks in nature and society are inherently directed in nature. These include the metabolic networks \cite{Jeong2000,Samal2006,Samal2011}, gene regulatory networks \cite{Milo2002}, signaling networks \cite{Maayan2005}, neural networks \cite{Sporns2004}, the world wide web (WWW) \cite{Broder2000}, online social networks \cite{Mcauley2012} and transportation networks \cite{Opsahl2010}. To this end, the two different discretizations of the Ricci curvature, Ollivier-Ricci and Forman-Ricci, have been developed and employed for the analysis of undirected complex networks to date \cite{Lin2010,Lin2011,Bauer2012,Jost2014,Loisel2014,Ni2015,Sandhu2015,Sandhu2016,Sreejith2016,Samal2018}. In a recent contribution \cite{Samal2018}, we have compared the performance of these two notions of Ricci curvature in several model and real-world undirected complex networks. To enable proper investigation of both model and real-world directed complex networks, we extend herein the notions of Forman-Ricci and Ollivier-Ricci curvatures to the domain of directed graphs or networks, and thereafter employ our definitions of the Forman-Ricci and Ollivier-Ricci curvature for a directed edge to analyze several model and real-world directed complex networks. Note that in an unpublished preprint \cite{Sreejith2016directed}, we had extended the definition of the Forman-Ricci curvature for an edge in undirected graphs to directed graphs, and in this contribution, we further broaden the notion of the Forman-Ricci curvature for a directed edge to also define an \textit{augmented Forman-Ricci curvature} for a directed edge which also takes into account the 2-dimensional faces arising from higher-order correlations between vertices in directed networks. We had also earlier developed the idea and explored the augmented Forman-Ricci curvature which accounts for two-dimensional faces in undirected complex networks \cite{Samal2018,Weber2017}. In addition, we also introduce here a fitting notion of the Ollivier-Ricci curvature for directed networks. Thus, the three notions of discrete Ricci curvature can hereafter be employed to investigate both undirected and directed complex networks.

The remainder of this paper is organized as follows. In the next section, we extend the two notions of Ricci curvature for undirected networks to directed networks. More precisely, we not only adapt the Ollivier-Ricci and Forman-Ricci curvatures for undirected networks discussed in Ref. \cite{Samal2018}, but we have also developed a directed analogue of the augmented Forman-Ricci curvature which also accounts for two- (and higher-) dimensional complexes. In the subsequent section, we list the datasets of directed networks, both model and real-world, employed to investigate the Forman-Ricci and Ollivier-Ricci curvature for directed networks. In the penultimate section, we describe our main results, and in the last section, we conclude with a summary and future outlook.


\section{Theoretical background}
\label{theory}

In this section, we extend the two discretizations of Ricci curvature, namely, Forman-Ricci and Ollivier-Ricci, which have been previously employed to study undirected networks to the setting of directed networks.

\subsection{Extension of Forman-Ricci curvature to directed networks}
\label{FormanCurvature}

Previously, we had adapted the Forman's discretization of Ricci curvature \cite{Forman2003} to the setting of undirected networks and have employed the curvature measure to investigate several model and real-world networks \cite{Sreejith2016,Sreejith2017}. The classical Ricci curvature operates directionally along vectors, and in the discrete setting of networks, Forman-Ricci curvature is associated with the discrete analogue of vectors, namely, edges in networks. Thus, Forman-Ricci curvature is an edge-based measure for geometrical characterization of networks. In an undirected network $G$, the Forman-Ricci curvature for an edge $e$ is given by the following formula \cite{Sreejith2016}:
\begin{equation}
\label{UndirectedFormanRicciEdge}
\mathbf{F}(e) = w_e \left( \frac{w_{v_1}}{w_e} +  \frac{w_{v_2}}{w_e}  - \sum_{e_{v_1},e_{v_2}\ \sim\ e} \left[\frac{w_{v_1}}{\sqrt{w_e w_{e_{v_1} }}} + \frac{w_{v_2}}{\sqrt{w_e w_{e_{v_2} }}} \right] \right)\,
\end{equation}
In Eq. \ref{UndirectedFormanRicciEdge}, $e$ denotes the edge under consideration between two vertices $v_1$ and $v_2$, $w_e$ denotes the weight of the edge $e$, and $w_{v_1}$ and $w_{v_2}$ denote the weights associated with the vertices $v_1$ and $v_2$, respectively. In Eq. \ref{UndirectedFormanRicciEdge}, $e_{v_1} \sim e$ and $e_{v_2} \sim e$ denote the set of edges incident on vertices $v_1$ and $v_2$, respectively, after \textit{excluding} the edge $e$ under consideration which connects the two vertices $v_1$ and $v_2$. Note that $e_{v_1},e_{v_2} \sim e$ under the summation symbol in Eq. \ref{UndirectedFormanRicciEdge} do not denote a double summation, but rather this notation enables compact representation of a single summation, and thus:
\begin{equation}
\label{UndirectedFormanRicciEdgeSingleSum}	
\sum_{e_{v_1},e_{v_2}\ \sim\ e} \left[\frac{w_{v_1}}{\sqrt{w_e w_{e_{v_1} }}} + \frac{w_{v_2}}{\sqrt{w_e w_{e_{v_2} }}} \right]= \sum_{e_{v_1}\ \sim\ e} \frac{w_{v_1}}{\sqrt{w_e w_{e_{v_1}}}} + \sum_{e_{v_2}\ \sim\ e} \frac{w_{v_2}}{\sqrt{w_e w_{e_{v_2} }}}
\end{equation}
From a geometrical perspective, Forman-Ricci curvature is a concept inspired from Riemannian and polyhedral geometry which quantifies for an edge $e$ the extent to which the network spreads out at the ends of edges in a network. The more is this spreading in the network, the more negative is the Forman-Ricci curvature of the network. Thus, Forman-Ricci curvature represents the information flow in the network along the edge $e$ and this physical interpretation can be easily seen from the defining formula (Eq. \ref{UndirectedFormanRicciEdge}). Moreover, using the definition of the Forman-Ricci curvature for an edge in the undirected network (Eq. \ref{UndirectedFormanRicciEdge}), one can elegantly define the Forman-Ricci curvature for a vertex $v$ in the undirected network as follows \cite{Sreejith2017}:
\begin{equation}
\label{UndirectedFormanRicciVertex}
\mathbf{F}(v) = \sum_{e \in E_{v}} \mathbf{F}(e) \,
\end{equation}
where $E_v$ denotes the set of edges incident on the vertex $v$. We remark that the Forman-Ricci curvature of a vertex represents a discrete analogue of scalar curvature in networks \cite{Sreejith2017}.

In the original work \cite{Forman2003}, Forman had not envisaged his curvature function for directed networks, but rather concentrated on more general geometric objects, the so called CW weighted cell complexes, of which polygonal meshes, for instance, represent one well known and widely employed example. Notice that in Forman's original work only positive weights were considered, since these weights represent therein generalizations of such natural geometric notions as length, area and volume. In consequence, no directed surfaces (or, more general complexes) can be modeled directly using Forman's original formalism. However, we will show here that it is simple to adapt the Forman-Ricci curvature for undirected networks \cite{Sreejith2016} to directed networks. Thus, this work, though restricted solely to networks, also represents a novel theoretical extension.

\begin{figure}
\centering
\includegraphics[width=.79\columnwidth]{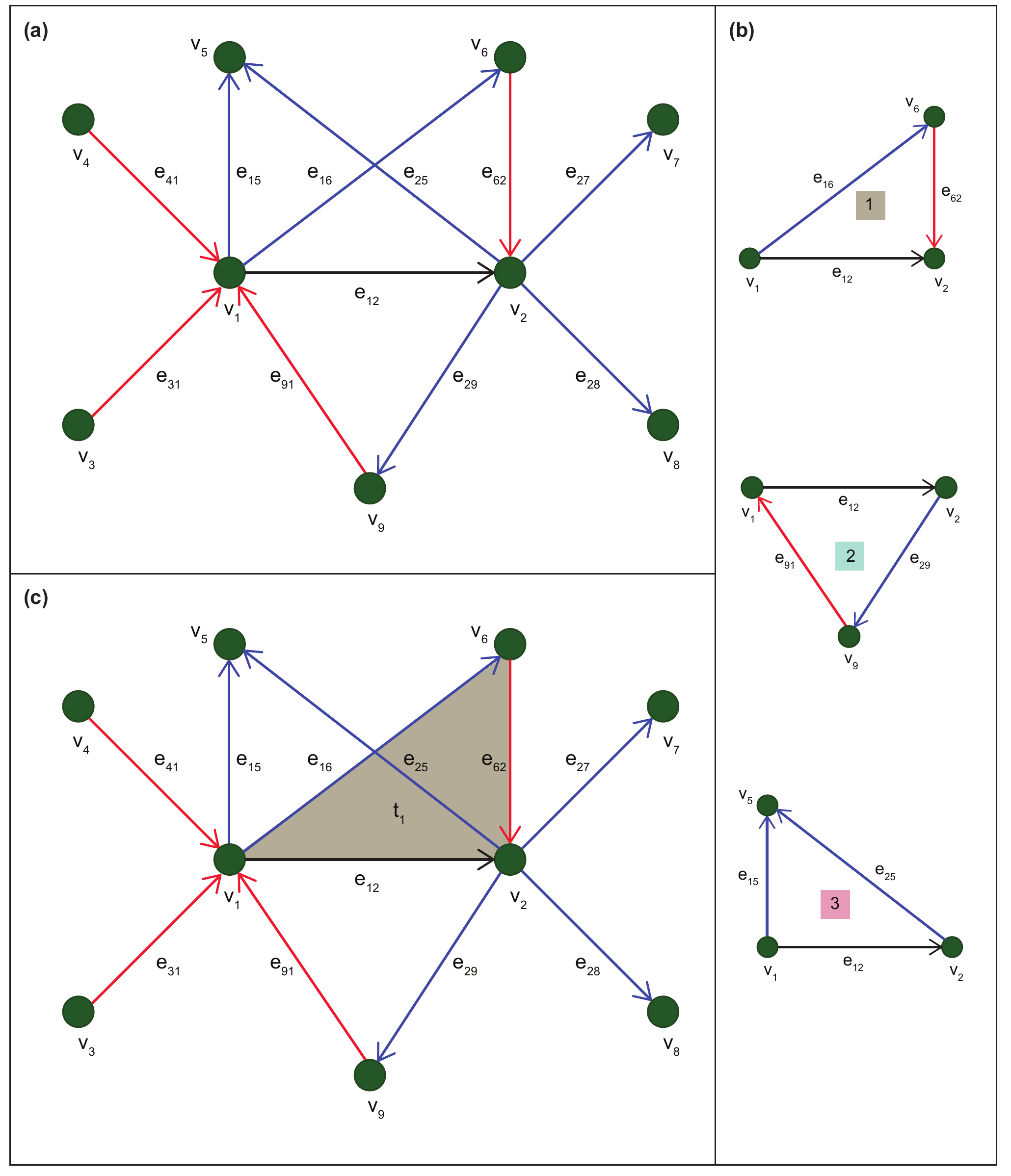}
\caption{Example network depicting the computation of the Forman-Ricci and the augmented Forman-Ricci curvature for directed edges. In this example network, the \textit{incoming} and \textit{outgoing} edges to vertices anchoring the edge $e_{12}$ under consideration (after excluding the edge $e_{12}$ itself) are shown in red and blue, respectively. (a) While computing the Forman-Ricci curvature for the directed edge $e_{12}$ which originates at vertex $v_1$ and terminates at vertex $v_2$, we consider only the incoming edges $e_{31}$, $e_{41}$ and $e_{91}$ at vertex $v_1$ and only the outgoing edges $e_{27}$, $e_{28}$ and $e_{29}$ at vertex $v_2$. (b) Shown here are the three essentially different types of directed triangles which can arise in directed networks. We here choose the first configuration which corresponds to a feed forward loop (FFL) motif as the directed triangular face or directed simplicial complex while computing the augmented Forman-Ricci curvature for a directed edge. (c) While computing the augmented Forman-Ricci curvature for the directed edge $e_{12}$ which originates at vertex $v_1$ and terminates at vertex $v_2$, we also consider the positive contribution from the directed triangular face $t_1$ formed by edges $e_{16}$, $e_{62}$ and $e_{12}$ in Eq. \ref{DirectedAugmentedFormanRicciEdge}, in addition to the negative contributions from the incoming edges $e_{31}$, $e_{41}$ and $e_{91}$ at vertex $v_1$ and the outgoing edges $e_{27}$, $e_{28}$ and $e_{29}$ at vertex $v_2$.}
\label{DirectedExplain}
\end{figure}

Notice that we can rearrange the terms in the definition of the Forman-Ricci curvature for an edge $e$ in Eq. \ref{UndirectedFormanRicciEdge} to separate the contributions of the two vertices $v_1$ and $v_2$ which anchor the edge $e$ under consideration as follows:
\begin{equation}
\label{UndirectedFormanRicciEdgeAlternate}
\mathbf{F}(e) = w_e \left( \frac{w_{v_1}}{w_e}   - \sum_{e_{v_1}\ \sim\ e} \frac{w_{v_1}}{\sqrt{w_e w_{e_{v_1} }}} \right)  + w_e \left( \frac{w_{v_2}}{w_e}   - \sum_{\ e_{v_2}\ \sim\ e} \frac{w_{v_2}}{\sqrt{w_e w_{e_{v_2} }}}  \right)\,
\end{equation}
Using the above equation, we define the Forman-Ricci curvature for a directed edge $e$ which originates from vertex $v_1$ and terminates at vertex $v_2$ in a directed network as follows:
\begin{equation}
\label{DirectedFormanRicciEdge}
\mathbf{F}(e) = w_e \left( \frac{w_{v_1}}{w_e}   - \sum_{e_{I,v_1}\ \sim\ e} \frac{w_{v_1}}{\sqrt{w_e w_{e_{I,v_1} }}} \right)  + w_e \left( \frac{w_{v_2}}{w_e}   - \sum_{\ e_{O,v_2} \ \sim\ e} \frac{w_{v_2}}{\sqrt{w_e w_{e_{O,v_2} }}}  \right)\,
\end{equation}
In Eq. \ref{DirectedFormanRicciEdge}, $e_{I,v_1} \sim e$ denote the set of incoming edges incident on vertex $v_1$ and $e_{O,v_2} \sim e$ denote the set of outgoing edges emanating from vertex $v_2$. Thus, while computing the Forman-Ricci curvature for a directed edge $e$ which originates from vertex $v_1$ and terminates at vertex $v_2$, we take into account only those directed edges that either terminate at vertex $v_1$ or originate at vertex $v_2$ (Fig. \ref{DirectedExplain}(a)). Moreover, while computing the Forman-Ricci curvature for a directed edge $e$, we ignore self-edges or self-loops at vertices $v_1$ and $v_2$ which can arise in real-world networks. In short, while computing the Forman-Ricci curvature for an edge $e$ in directed networks, we consider only those neighboring edges whose direction is compatible with the direction of edge $e$, and moreover, this fundamentally ensures that the resulting flow is consistent with the direction of edge $e$ under consideration (Fig. \ref{DirectedExplain}(a)). We remark that our definition of the Forman-Ricci curvature for a directed edge is, by no means, the only possible choice. Beyond the mathematical reason discussed above, it is the most natural definition for the modelling of diverse real networks considered here.

For a given vertex $v$ in a directed network, one can distinguish between its incoming and outgoing edges. Given a vertex $v$, let us denote the set of \textit{incoming} and \textit{outgoing} edges for a vertex $v$ by $E_{I,v}$ and $E_{O,v}$, respectively (Fig. \ref{DirectedExplain}(a)). Then, one can elegantly define the \textit{In Forman-Ricci curvature} $\mathbf{F}_I(v)$ and the \textit{Out Forman-Ricci curvature} $\mathbf{F}_O(v)$ as follows:
\begin{equation}
\label{DirectedFormanRicciVertexIn}
\mathbf{F}_I(v) = {\sum_{e \in E_{I,v}} \mathbf{F}(e)}\,;
\end{equation}
\begin{equation}
\label{DirectedFormanRicciVertexOut}
\mathbf{F}_O(v) = {\sum_{e \in E_{O,v}} \mathbf{F}(e)}\,;
\end{equation}
where the summations are taken over only the incoming and outgoing edges, respectively. Moreover, one can obtain the total amount of flow through a vertex $v$ in the directed network as follows:
\begin{equation}
\label{DirectedFormanRicciVertexFlow}
\mathbf{F}_{T}(v) =  \mathbf{F}_I(v) - \mathbf{F}_O(v)\,.
\end{equation}
Thus, in a directed network, one can associate three Forman-Ricci curvatures with each vertex $v$, namely, In Forman-Ricci curvature $\mathbf{F}_I(v)$, the Out Forman-Ricci curvature $\mathbf{F}_O(v)$ and the total flow $\mathbf{F}_{T}(v)$.

\subsection{Augmented Forman-Ricci curvature for directed networks}
\label{AugmentedFormanCurvature}

Traditional network measures consider only the pairwise correlations between vertices as embodied by edges in a complex network while ignoring higher-order correlations between $n$-tuples of vertices. Such higher-order correlations naturally arise between vertices in real-world networks such as scientific co-authorship networks \cite{Newman2004} and disease-gene association networks \cite{Goh2007}. From a geometric point of view, these higher-order correlations between vertices in a network are captured by two-dimenesional faces determined by the relevant $n$-tuples of vertices. For instance, one may construct two-dimensional polyhedral complexes by inserting a triangle into any connected triple of vertices or cycle of length 3, a quadrangle into a cycle of length 4, a pentagon into a cycle of length 5, and so on. Notably, the original definition \cite{Forman2003} of the Forman-Ricci curvature for general $CW$ complexes applies to complexes of dimension higher or equal to 2, which is given by:
\begin{equation}
\label{FormanRicciComplex}
{\rm F}^{\#} (e) = w_e \left[ \sum_{e < f} \frac{w_e}{w_f}+\sum_{v < e} \frac{w_v}{w_e} \right. - \left. \sum_{\hat{e} \parallel e} \left| \sum_{\hat{e},e < f} \frac{\sqrt{w_e \cdot w_{\hat{e}}}}{w_f} - \sum_{v 	 < e, v < \hat{e}} \frac{w_v}{\sqrt{w_e \cdot w_{\hat{e}}}} \right| \right] \; ;
\end{equation}
where $w_e$ denotes the weight of edge $e$, $w_v$ denotes the weight of vertex $v$, $w_f$ denotes the weight of face $f$, $\sigma < \tau$ means that $\sigma$ is a face of $\tau$, and $||$ signifies \textit{parallelism}, i.e. the two cells have a common \textit{parent} (higher dimensional face) or a common \textit{child} (lower dimensional face), but not both a common parent and common child. In particular, the above formula (Eq. \ref{FormanRicciComplex}) for complexes having solely triangular faces $t$ while neglecting correlations between more than 3 vertices becomes:
\begin{equation}
\label{UndirectedAugmentedFormanRicciEdge}
{\rm F}^{\#} (e) = w_e \left[ \sum_{e < t} \frac{w_e}{w_t} + (\frac{w_{v_1}}{w_e} + \frac{w_{v_2}}{w_e}) \right. - \left. \sum_{e_{v_1},e_{v_2}\ \sim\ e,  e_{v_1},e_{v_2}\ \nless\ t} \left( \frac{w_{v_1}}{\sqrt{w_e \cdot w_{e_{v_1}}}} + \frac{w_{v_2}}{\sqrt{w_e \cdot w_{e_{v_2}}}} \right) \right] \;
\end{equation}
We denote the resulting quantity given by the above equation as the \textit{augmented Forman-Ricci curvature} for an edge $e$ which also accounts for two-dimensional simplicial complexes of length 3 arising in undirected networks \cite{Samal2018}. We remark that in a unweighted and undirected network $G$ with $w_t = w_e = w_v = 1, \; \forall t \in T(G), e \in E(G), v \in V(G)$, where $T(G)$, $E(G)$ and $V(G)$ represent the set of triangular faces, edges and vertices, respectively, there is a simple relationship \cite{Samal2018} between Forman-Ricci curvature ${\mathrm F} (e)$ and augmented Forman-Ricci curvature ${\mathrm F}^{\#} (e)$  for an edge $e$ which is:
\begin{equation}
\label{UndirectedUnweightedAugmentedFormanRicciEdge}
{\rm F}^{\#} (e) = {\rm F} (e) + 3m
\end{equation}
where $m$ is the number of triangles containing edge $e$ under consideration in the undirected network. Importantly, the original definition of the Forman-Ricci curvature concerns $CW$ complexes of dimension at least 2 (Eq. \ref{FormanRicciComplex}), and thus, our definition of the augmented Forman-Ricci curvature naturally extends to hypergraphs. In a future contribution, we will report ongoing work on application of the augmented Forman-Ricci curvature (Eq. \ref{UndirectedUnweightedAugmentedFormanRicciEdge}) to complex hypernetworks.

We now extend the definition of the augmented Forman-Ricci curvature for an edge $e$ in undirected networks to directed networks. This adaptation is rather straightforward provided once the appropriate choice for directed triangular faces is made. While making the necessary choice, it is imperative to remark that there are three essentially different types of directed triangles which can arise in complex networks (Fig. \ref{DirectedExplain}(b)), and of these three possible choices, we opt for the first configuration to be the directed simplicial complex following the recent work by Masulli \& Villa \cite{Masulli2016} and Courtney \& Bianconi \cite{Courtney2018}. Note that the first configuration shown in Fig. \ref{DirectedExplain}(b) and chosen as the directed simplicial complex is the extensively studied feed forward loop (FFL) motif which is preponderant in several real-world directed networks including transcriptional regulatory networks \cite{Milo2002}. We highlight that our choice of the faces of the directed simplicial complex is by no means universal, and in particular, from a purely geometric or topological point of view, the second configuration in Fig. \ref{DirectedExplain}(b) corresponding to a cycle could also represent a viable choice. However, we followed Masulli \& Villa \cite{Masulli2016} and Courtney \& Bianconi \cite{Courtney2018} in choosing the first configuration in Fig. \ref{DirectedExplain}(b) as the faces of the directed simplicial complex because this configuration is in concordance with the direction of the information flow along the edge $e$ under consideration, and thus, conforms with the intuition of the Forman-Ricci curvature as a measure of geodesic dispersal and flow through an edge. With the above-mentioned choice of the face of the directed simplicial complex (Fig. \ref{DirectedExplain}(c)), the mathematical formula for the augmented Forman-Ricci curvature for a directed edge $e$ which originates from vertex $v_1$ and terminates at vertex $v_2$ in a directed network is given by:
\begin{equation}
\label{DirectedAugmentedFormanRicciEdge}
{\rm F}^{\#} (e) = w_e \left[ \left( \sum_{e < t} \frac{w_e}{w_t} + (\frac{w_{v_1}}{w_e} + \frac{w_{v_2}}{w_e}) \right) \right. - \left. \sum_{e_{I,v_1},e_{O,v_2}\ \sim\ e,  e_{I,v_1},e_{O,v_2}\ \nless\ t} \left( \frac{w_{v_1}}{\sqrt{w_e \cdot w_{e_{I,v_1}}}} + \frac{w_{v_2}}{\sqrt{w_e \cdot w_{e_{O,v_2}}}} \right) \right] \;
\end{equation}
In this contribution, we have also explored the augmented Forman-Ricci curvature for a directed edge $e$ in model and real-world directed networks.


\subsection{Extension of Ollivier-Ricci curvature to directed networks}
\label{OllivierCurvature}

We here briefly describe the Ollivier-Ricci curvature and its adaptation to directed networks. Ollivier's discretization of Ricci curvature stems from the essential observation that ``In positive curvature, balls are closer than their centers are, while in negative curvature, balls are farther than their centers are'' \cite{Ollivier2013}. Indeed, in a $n$-dimensional Riemannian manifold, this property is captured precisely by the Ricci curvature ${\rm Ric}(v)$ in the tangential direction $v$, as quantized in the following formula for the average distance $\delta = d(x,y)$ between two infintesimal balls $B_x$ and $B_y$ with centers $x$ and $y$, respectively, and radii equal to $\varepsilon$:
\begin{equation}
\label{OllivierRicci-1}
\delta\left(1 - \frac{\varepsilon^2}{2(n+2)}{\rm Ric}(v) + O(\varepsilon^3 + \varepsilon^2\delta) \right)\,,
\end{equation}
where $\varepsilon, \delta \rightarrow 0$. It is most natural to extend the above definition from volumes of infintesimal balls to more general measures. The additional ingredient needed to this end is to find a good definition for the distance between two measures. The most natural as well as easy to employ such distance is the Monge-Kantorovich-Wasserstein transportation metric $W_1$ (see, e.g. \cite{Lott2009}) which for graphs or networks has the following form, $W_1(m_x, m_y)$, i.e. the transportation distance between the two probability measures $m_x$ and $m_y$, is given by:
\begin{equation}
\label{OllivierRicci-2}
W_1(m_x, m_y)=\inf_{\mu_{x,y}\in \prod(m_x, m_y)}\sum_{(x',y')\in V\times V}d(x', y')\mu_{x,y}(x', y'),
\end{equation}
where $m_x$, $m_y$ represent two probability measures associated with the vertices $x$, $y$, respectively, and where $\prod(m_x, m_y)$ denotes the set of probability measures $\mu_{x,y}$ that satisfy:
\begin{equation}
\label{OllivierRicci-3}
\sum_{y'\in V}\mu_{x,y}(x', y')=m_x(x'), \,\,\sum_{x'\in V}\mu_{x,y}(x', y')=m_y(y').
\end{equation}
Note that the measure $m$, associated to the set of vertices of a network is, evidently, a discrete measure.

Then Ollivier-Ricci curvature \cite{Ollivier2007,Ollivier2009} is defined as follows:
\begin{equation}
\label{OllivierRicci-4}
\kappa(x,y) = 1 - \frac{W_1(m_x,m_y)}{d(x,y)}\,
\end{equation}
where $m_x, m_y$ represent the measures of the balls around $x$ and $y$, respectively.

We highlight that while computing the Wasserstein distance \cite{Vaserstein1969} or earth mover's distance $W_1(m_x, m_y)$ for an edge $e=(x,y)$ in an undirected and unweighted network, the unit probability measure $m_x$ concentrated at the vertex $x$ is equally distributed among the immediate neighbors of $x$ while the unit probability measure $m_y$ concentrated at the vertex $y$ should also be equally distributed among the immediate neighbors of $y$, and this is used to define the transportation plan from $x$ to $y$ \cite{Ni2015}. Moreover, the ground distance $d(x,y)$ for an edge $e=(x,y)$ in an undirected and unweighted network is equal to the shortest path between vertices $x$ and $y$ which is trivially equal to 1 \cite{Ni2015}. Thus, the Ollivier-Ricci curvature for an edge $e=(x,y)$ in an undirected and unweighted network is given by:
\begin{equation}
\label{UndirectedOllivierRicciEdge}
\kappa(x,y) = 1 - W_1(m_x,m_y).
\end{equation}

In order to extend the notion of Ollivier-Ricci curvature to directed networks, we compute the analogue of Wasserstein distance $\overline{W_1}(m_x, m_y)$ for a directed edge $\bar{e}=(\overline{x,y})$ in a directed and unweighted network as follows. The unit probability measure $m_x$ concentrated at the starting vertex $x$ is equally distributed among the immediate neighbours with incoming edges to $x$ while the unit probability measure $m_y$ concentrated at the terminating vertex $y$ is equally distributed among the immediate neighbors with outgoing edges from $y$. In other words, similar to the definition of the Forman-Ricci curvature for a directed edge $\bar{e}$, while computing the Ollivier-Ricci curvature, we consider only the incoming edges to the starting vertex $x$ and only the outgoing edges from the terminating vertex $y$. Thus, the Ollivier-Ricci curvature for a directed edge $\bar{e}=(\overline{x,y})$ in a directed and unweighted network is given by:
\begin{equation}
\label{DirectedOllivierRicciEdge}
\kappa(\overline{x,y}) = 1 - \overline{W_1}(m_x,m_y).
\end{equation}
We remark that the directed analogue of Wasserstein distance is no longer a proper metric as the $\overline{W_1}(m_y, m_x)$ may not even be defined in the directed network.

Similar to the case of Forman-Ricci curvature in directed networks, one can elegantly define the \textit{In Ollivier-Ricci curvature} $\mathbf{O}_I(v)$ and the \textit{Out Ollivier-Ricci curvature} $\mathbf{O}_O(v)$ as follows:
\begin{equation}
\label{DirectedOllivierRicciVertexIn}
\mathbf{O}_I(v) = {\sum_{e \in E_{I,v}} \kappa(e)}\,;
\end{equation}
\begin{equation}
\label{DirectedOllivierRicciVertexOut}
\mathbf{O}_O(v) = {\sum_{e \in E_{O,v}} \kappa(e)}\,;
\end{equation}
where the summations are taken over only the incoming and outgoing edges, respectively. Moreover, one can obtain the total amount of flow through a vertex $v$ in the directed network as follows:
\begin{equation}
\label{DirectedOllivierRicciVertexFlow}
\mathbf{O}_{T}(v) =  \mathbf{O}_I(v) - \mathbf{O}_O(v)\,.
\end{equation}

Before concluding this section, let us note that the Ollivier-Ricci curvature has no immediate generalization to hypernetworks \cite{Banerjee2017,Asoodeh2018}. However, one proposition to compute the Ollivier-Ricci curvature in hypernetworks is as follows. One can associate to any two-dimensional cell complex with its dual graph, with vertex weights equal to the original face weights and with the edge weights coinciding to those in the original complex. Finally, one can compute the Ollivier-Ricci curvature of the hypernetwork by computing the Ollivier-Ricci curvature of the associated dual graph.


\section{Datasets}
\label{datasets}

We have analyzed here the definitions introduced in the previous section of the Forman-Ricci, the Augmented Forman-Ricci and the Ollivier-Ricci curvature for a directed edge, in model and real-world directed networks.

We have considered two generative models of directed networks, namely, Erd\"{o}s-R\'{e}nyi (ER) model \cite{Erdos1961} and Scale-free (SF) model from vertex fitness scores \cite{Chung2002,Goh2001}. ER model generates an ensemble $G(n,p)$ of random networks where $n$ is the number of nodes and $p$ is the probability that each possible directed edge exists between any pair of vertices in the network. SF model from vertex fitness scores generates directed networks with power-law distribution for both in-degree and out-degree of vertices based on four parameters, namely, $n$, $m$, $\lambda_{in}$ and $\lambda_{out}$, where $n$ is the number of vertices, $m$ is the number of edges, $\lambda_{in}$ is the exponent of the in-degree distribution of vertices and $\lambda_{out}$ is the exponent of the out-degree distribution of vertices in the network.

In addition, we have also considered the following directed real-world networks. Air traffic control \cite{Kunegis2013} is a directed network which was reconstructed based on the Preferred Routes Database of the US Federal Aviation Administration National Flight Data Center (NFDC). Air traffic control network contains 1226 vertices and 2613 directed edges where vertices correspond to airports or service centers and directed edges represent preferred routes recommended by the NFDC. Twitter List \cite{Mcauley2012} is a directed network of connections between Twitter users which contains 23370 vertices and 33101 directed edges where vertices correspond to users and directed edges signify that one user follows another user on Twitter. \textit{E. coli} TRN \cite{Salgado2013} is the transcriptional regulatory network in the bacterium \textit{Escherichia coli} which contains 3073 vertices and 7853 directed edges where vertices corresponds to genes and directed edges represent control of target gene expression through transcription factors. \textit{B. subtilis} TRN \cite{Kumar2015} is the transcriptional regulatory network in bacterium \textit{Bacillus subtilis} which contains 1594 vertices and 2902 directed edges where vertices correspond to genes and directed edges represent control of target gene expression through transcription factors. Human protein interaction network \cite{Ewing2007} is a directed network of protein-protein interactions in human cells which contains 2239 vertices and 6452 directed edges where vertices corresponds to proteins and directed edges represent interactions between two proteins. Phosphonetwork \cite{Hu2013} is a human protein phosphorylation network created using combined bioinformatics and protein microarray based strategy. Phosphonetwork contains 1290 vertices and 4372 directed edges where vertices correspond to human protein kinases or the substrate proteins of human protein kinases and directed edges represent phosphorylation of the substrate protein by the human protein kinases.

Supplementary Table S1 lists the different model and real-world directed networks analyzed here. In addition, Supplementary Table S1 lists several network measures such as maximum degree, average degree, average in-degree, average out-degree and size of the largest weakly connected component for the different model and real-world directed networks analyzed here. Furthermore, we have chosen different combinations of input parameters for the ER and SF models considered here to generate directed networks with different average degree. Also, we have sampled 100 directed networks starting with different random seed for a specific combination of input parameters from the two models of directed networks, and the results reported here for model networks is an average over the sample of 100 directed networks with chosen input parameters. Note that we have considered here directed yet unweighted model and real-world networks, and thus, the weights of vertices, edges and faces are taken to be 1 while computing the discrete Ricci curvatures.


\section{Results and discussion}
\label{results}

\subsection{Edge curvature in directed networks}

\begin{figure}
\centering
\includegraphics[width=.97\columnwidth]{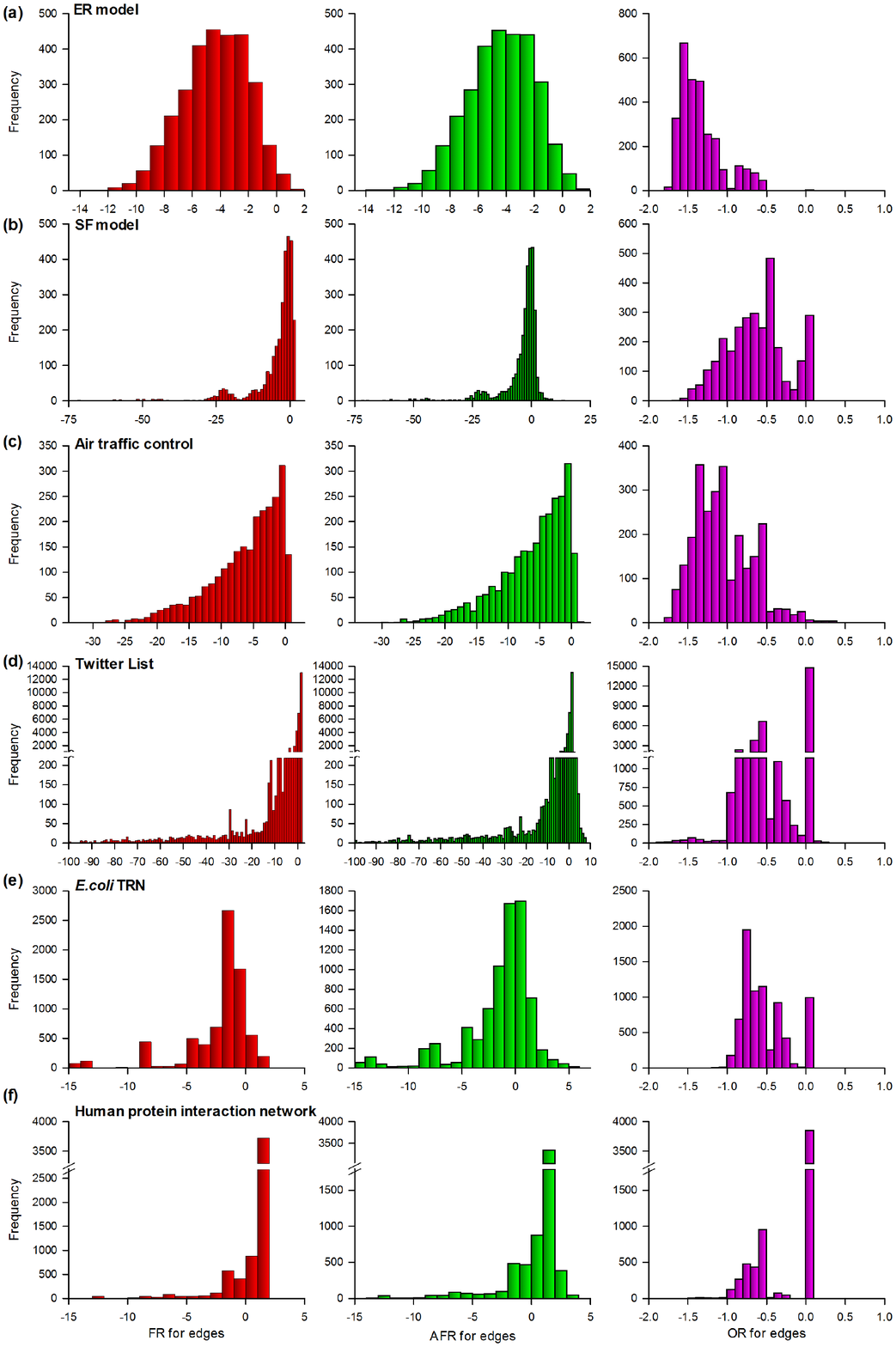}
\caption{Distribution of the Forman-Ricci (FR), the Augmented Forman-Ricci (AFR) and the Ollivier-Ricci (OR) curvature for directed edges in model and real-world directed networks. (a) ER model with $n=1000$, $p=0.003$. (b) SF model with $n=1000$, $m=3000$, $\lambda_{in}$=2.1, $\lambda_{out}$=2.1. (c) Air traffic control. (d) Twitter List. (e)  \textit{E. coli} TRN. (f) Human protein interaction network. }
\label{dist_edge}
\end{figure}

In Fig. \ref{dist_edge}, we show the distribution of the Forman-Ricci (FR), the Augmented Forman-Ricci (AFR) and the Ollivier-Ricci (OR) curvature for edges in model and real-world directed networks. Note that the OR curvature for edges in directed and unweighted networks takes value in the range -2 to +1 which is independent of the number of edges or vertices in the network. In contrast, the FR or AFR curvature for edges depend on the number of edges or vertices in the directed network. Specifically, the FR or AFR curvature for an edge in directed and unweighted networks can take a very high negative value if the in-degree of the originating vertex and the out-degree of the terminating vertex of the directed edge under consideration is very high, and this can especially happen in large directed complex networks. From Fig. \ref{dist_edge}, it is seen that, in both model and real-world directed networks considered here, the FR or AFR or OR curvature for most edges is negative. Recall that Ricci curvature controls the growth of volumes in the classical Riemannian setting \cite{Berger2000,Jost2017}. Thus, spaces with negative Ricci curvature have exponential type growth while those with positive Ricci curvature have a finite diameter, and this result also holds for both the FR and OR curvature \cite{Forman2003,Ollivier2009}. Thus, our observation that most edges in considered real-world directed networks have a negative FR or AFR or OR curvature suggests that these networks have potential of infinite growth.

From Fig. \ref{dist_edge}, it is also seen that one can clearly distinguish the two models of directed networks, ER and SF, by the observed nature of the curvature distribution for edges. In directed ER networks, the distribution of the FR or AFR curvature for edges is narrow. But, in directed SF networks, the distribution of the FR or AFR curvature for edges is broad. In Fig. \ref{dist_edge}, we also show the distribution of the FR or AFR curvature for edges in four real directed networks, namely, Air traffic control, Twitter List, \textit{E. coli} TRN and Human protein interaction network. Similar to the directed SF networks, the four real directed networks have a broad distribution of the FR or AFR curvature for edges. These observations of broad distribution of the FR or AFR curvature for edges in real directed networks is consistent with the observed scale-free architecture with power-law degree distribution in considered networks. Also, note that more number of directed edges have positive AFR curvature in comparison to those with positive FR curvature in model and real directed networks considered here (Fig. \ref{dist_edge}), and this can be explained by positive contributions to the curvature by the directed simplicial complexes which are accounted in the definition of the AFR curvature. Although in directed and unweighted networks, the distribution of the OR curvature for edges is bounded by the same range, we still find that the distribution of the OR curvature for edges in directed ER networks is very different from that for directed SF networks (Fig. \ref{dist_edge}). Interestingly, the mean of the distribution of the OR curvature for edges in the real directed networks is in most cases found to be closer to that of the directed SF network than the directed ER network (Fig. \ref{dist_edge}).

\begin{figure}
\centering
\includegraphics[width=.97\columnwidth]{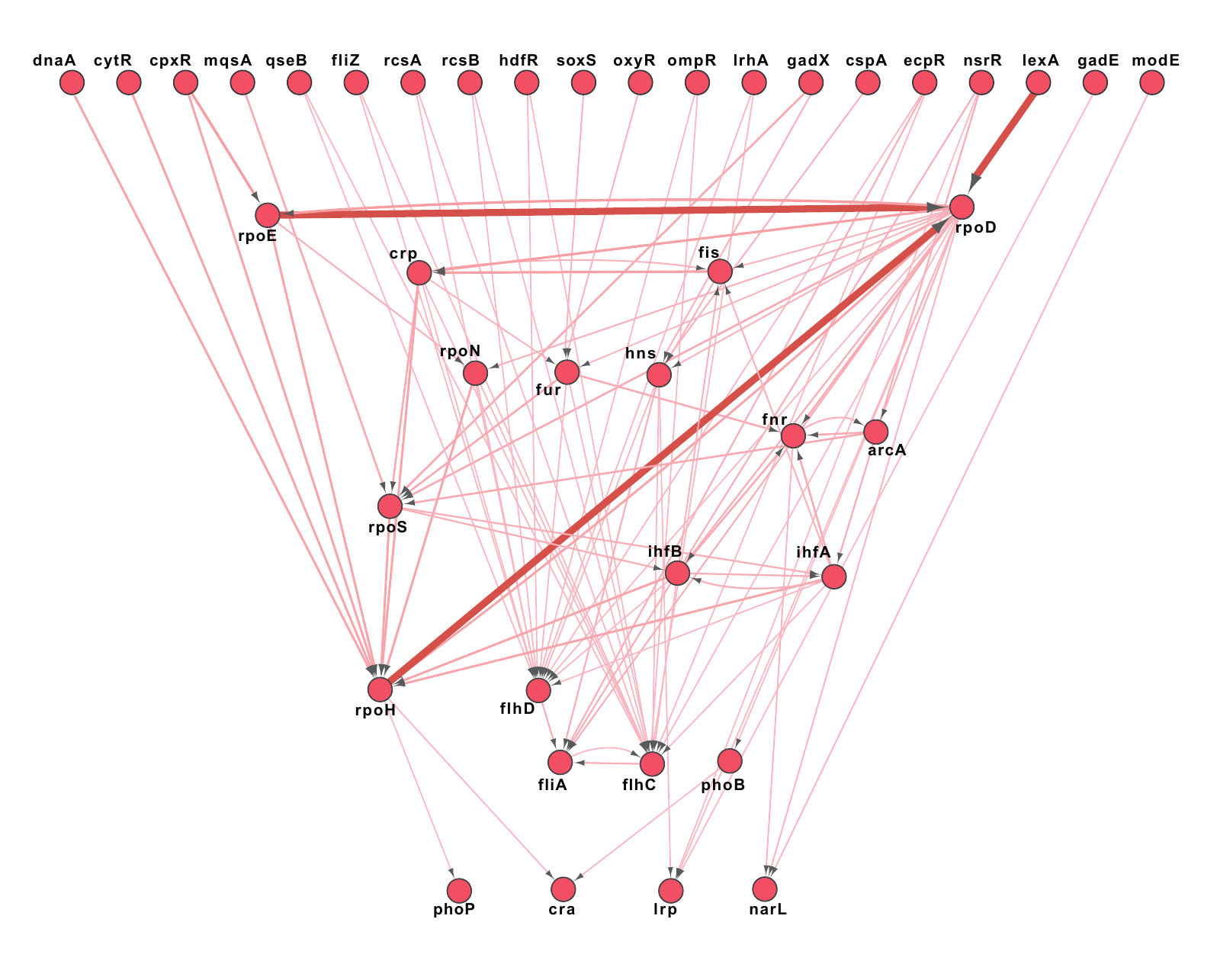}
\caption{Visualization of the subnetwork within the \textit{E. coli} transcriptional regulatory network (TRN) containing vertices connected by directed edges with highly negative augmented Forman-Ricci curvature (AFR), i.e., edges with ${\mathrm F}^{\#}(e)\le-58$. In this figure, the width of the directed edges is proportional to the magnitude of its AFR curvature. Interestingly, sigma factors including rpoD, rpoH and rpoE which play a critical role in transcriptional initiation are connected by directed edges with very high negative AFR curvture in this core subnetwork of the \textit{E. coli} TRN.}
\label{ecoli_afr}
\end{figure}

In Fig. \ref{ecoli_afr}, we show the portion of the \textit{E. coli} transcriptional regulatory network (TRN) which contains only vertices connected by directed edges with highly negative AFR curvature, i.e., edges with ${\mathrm F}^{\#}(e)\le-58$. Interestingly, sigma factors such as rpoD, rpoH and rpoE, which play a critical role in transcription initiation, are connected by directed edges with very high negative AFR curvature in this core subnetwork of the \textit{E. coli} TRN (Fig. \ref{ecoli_afr}).

We next compared the OR with FR and AFR curvature for edges in model and real-world directed networks (Table \ref{tab-1}). In directed ER and SF networks, we find a high positive correlation between the OR and FR or AFR curvature for edges when the model directed networks are sparse with small average degree, but the observed correlation vanishes with increase in average degree of the model directed networks (Table \ref{tab-1}; Supplementary Table S1). Note that the FR and AFR curvature for edges in directed ER and SF networks have very similar correlation with OR curvature (Table \ref{tab-1}). We have also compared the OR with FR and AFR curvature of edges in six real-world directed networks considered here (Table \ref{tab-1}). In majority of the real-world directed networks considered here, we find a high positive correlation between the OR and FR or AFR curvature for edges (Table \ref{tab-1}). Both the FR and AFR curvature for edges also have similar correlation with OR curvature in real-world directed networks considered here (Table \ref{tab-1}). Hence, we find a high positive correlation between the OR and FR or AFR curvature for edges in sparse model directed networks and several real-world directed networks. Note that, in a recent contribution \cite{Samal2018}, we have reported very high correlation between the OR and FR or AFR curvature for edges in sparse model undirected networks and several real-world undirected networks. Together, our results reinforce our preference towards Forman's discretization of the Ricci curvature for analysis of large complex networks. We highlight that Forman's discretization of Ricci curvature and its augmented version is very simple to compute in large undirected or directed complex networks, while the corresponding computation for Ollivier's discretization may be prohibitive in many large undirected or directed networks.

As emphasized above, the FR, AFR and OR curvatures are edge-centric measures for the analysis of complex networks \cite{Sreejith2016,Samal2018}. In network theory, edge-centric measures are less common than vertex-centric measures, and the present toolkit of other edge-centric measures for the analysis of undirected networks include edge betweenness centrality \cite{Freeman1977}, embeddedness \cite{Marsden1984} and dispersion \cite{Backstrom2014}. Among these edge-centric measures for undirected networks, edge betweenness centrality can also be used to analyze directed networks. Edge betweennness centrality \cite{Freeman1977,Girvan2002,Newman2010} is a global measure that quantifies the number of shortest paths that pass through an edge in an undirected or directed network. Edges with high betweenness centrality are bottlenecks for flows in the network. We have investigated the correlation between the edge-centric curvatures, OR, FR and AFR, with edge betweenness centrality in model and real directed networks (Table \ref{tab-2}). In directed ER and SF networks, we find a high negative correlation between edge betweenness centrality and OR or FR or AFR curvature for edges, especially, when the model directed networks are sparse with small average degree (Table \ref{tab-2}). Similarly, in majority of the real-world directed networks considered here, we find a moderate to high negative correlation between edge betweenness centrality and OR or FR or AFR curvature for edges (Table \ref{tab-2}). We remark that both FR and AFR curvature for a directed edge are local measures dependent on the weights of anchoring vertices and neighboring edges in the network. In contrast, the edge betweenness centrality is a global measure dependent on all shortest paths in the directed network. Thus, FR and AFR curvature for edges are simpler and faster to compute in large networks in comparison to edge betweenness centrality.


\subsection{Vertex curvature in directed networks}

Using the definition of the OR, FR or AFR curvature for directed edges, it is simple to define two curvatures associated with each vertex in a directed network (See section \ref{theory}). Thus, for a given vertex in a directed network, one can compute the In Ollivier-Ricci (IOR), Out Ollivier-Ricci (OOR), In Forman-Ricci (IFR), Out Forman-Ricci (OFR), In Augmented Forman-Ricci (IAFR) and Out Augmented Forman-Ricci (OAFR) curavture (See section \ref{theory}). We have also compared the IOR with IFR and IAFR curvature and OOR with OFR and OAFR curvature in model and real-world directed networks (Table \ref{tab-3}). In directed ER and SF networks, IFR and IAFR curvature have a high positive correlation with IOR curvature, OFR and OAFR curvature have a high positive correlation with OOR curvature, and moreover, the observed correlation is less sensitive to average degree of model directed networks unlike the observed correlation between OR and FR or AFR curvature for edges (Table \ref{tab-3}). Also, in majority of the real directed networks considered here, IFR and IAFR curvature have a moderate to high positive correlation with IOR curvature, and OFR and OAFR curvature have a moderate to high positive correlation with OOR curvature (Table \ref{tab-3}). Taken together, there is a moderate to high positive correlation between the two discretizations of Ricci curvature due to Forman and Ollivier, in most of the analyzed directed networks.

In-degree gives the number of edges incident to the vertex while out-degree gives the number of edges originating from the vertex in a directed network. We have also compared the in-degree with IOR, IFR and IAFR curvature and out-degree with OOR, OFR and OAFR curvature in model and real-world directed networks (Table \ref{tab-4}). In directed ER and SF networks, the in-degree has high negative correlation with IOR, IFR and IAFR curvature, and the out-degree has high negative correlation with OOR, OFR and OAFR curvature, and the observed correlation is less sensitive to average degree of the model directed networks. In majority of the real directed networks considered here, the in-degree and out-degree have moderate to high negative correlation with IOR and OOR curvature, respectively (Table \ref{tab-4}). In the real directed networks considered here, the out-degree has high negative correlation with OFR or OAFR curvature (Table \ref{tab-4}). However, in some of the real directed networks considered here, the in-degree has no or weak negative correlation with IFR or IAFR curvature (Table \ref{tab-4}).

Betweenness centrality \cite{Freeman1977} of a vertex gives the fraction of shortest paths between all pairs of vertices in the network that pass through that vertex. In directed ER and SF networks, betweenness centrality has moderate to high negative correlation with IOR, OOR, IFR, OFR, IAFR and OAFR curvature (Supplementary Table S2). Moreover, in majority of the real directed networks considered here, betweenness centrality has moderate to high negative correlation with IOR, OOR, IFR, OFR, IAFR and OAFR curvature (Supplementary Table S2). Pagerank \cite{Page1999,Langville2005} is an algorithm for directed networks which was originally developed to rank websites by the Google search engine. Pagerank can be employed with any directed network to measure importance of different vertices in the network. The algorithm counts the number and quality of incoming edges to a vertex to estimate the importance of a vertex in the network. In directed ER and SF networks, page rank has a high negative correlation with IOR, IFR and IAFR curvature, while page rank has no or weak negative correlation with OOR, OFR and OAFR curvature (Supplementary Table S3). In contrast, page rank has no or weak negative correlation with IOR, IFR, IAFR, OOR, OFR and OAFR curavture in majority of the real directed networks considered here (Supplementary Table S3).


\subsection{Curvature and robustness of directed networks}
\label{robustness}

\begin{figure}
\centering
\includegraphics[width=.97\columnwidth]{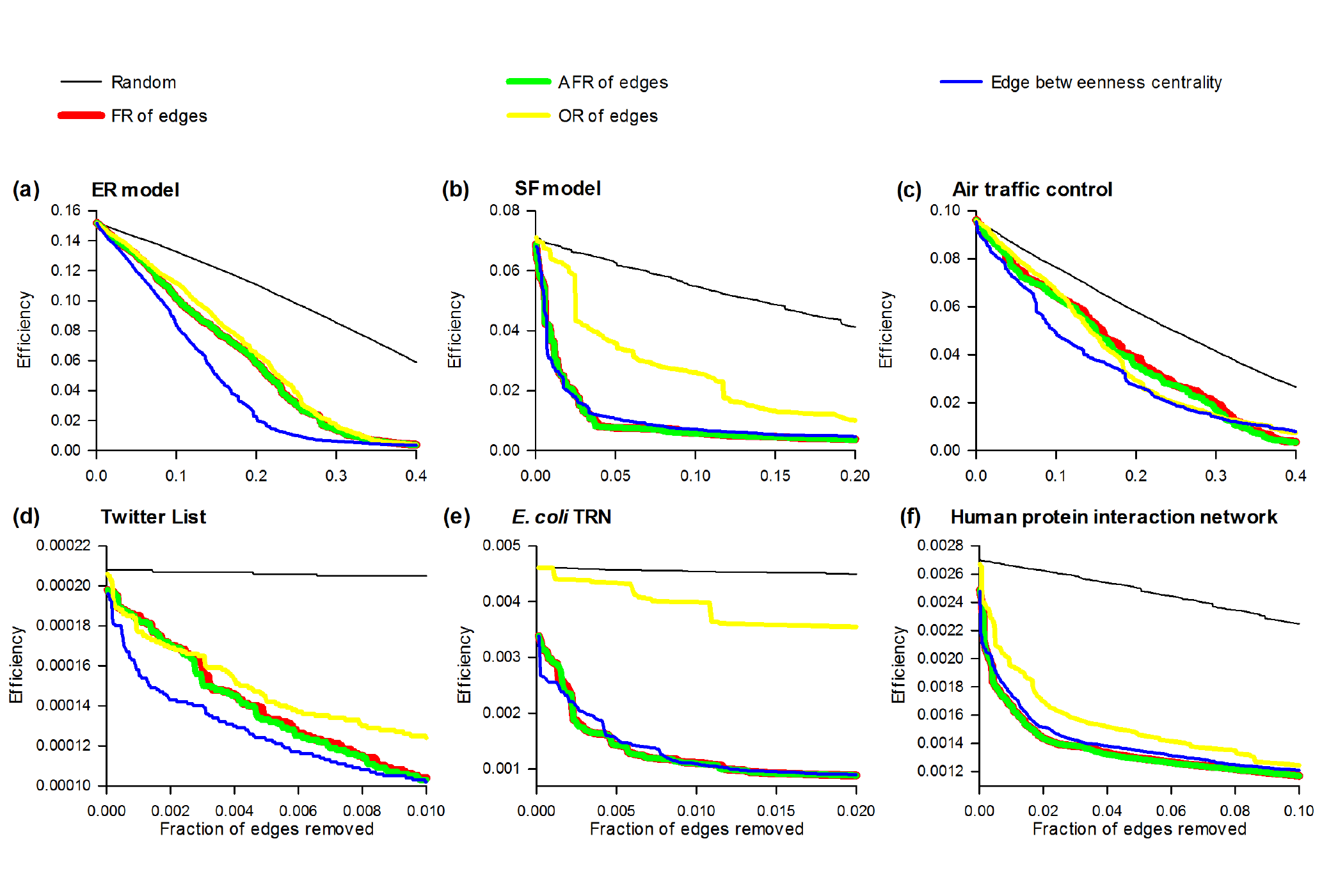}
\caption{Communication efficiency as a function of the fraction of edges removed in model and real-world directed networks. (a) ER model with $n=1000$, $p=0.003$. (b) SF model with $n=1000$, $m=3000$, $\lambda_{in}$=2.1, $\lambda_{out}$=2.1. (c) Air traffic control. (d) Twitter List. (e)  \textit{E. coli} TRN. (f) Human protein interaction network.}
\label{rob_edge}
\end{figure}

We next investigated the effect of removing directed edges based on increasing order of their OR, FR or AFR curvature on the large-scale connectivity of directed networks. Communication efficiency \cite{Latora2001} is a measure that captures how efficiently the information can be exchanged across the network, and the measure can be used to quantify a network's resistance to failure in face of small perturbations. Fig. \ref{rob_edge} shows the communication efficiency in model and real-world directed networks as a function of the fraction of edges removed. Here, the order of removing edges is based on the following criteria: (a) random order, (b) increasing order of OR curvature, (c) increasing order of FR curvature, (d) increasing order of AFR curvature, and (e) decreasing order of edge betweenness centrality. It is clear that targeted removal of edges with highly negative OR, FR or AFR curvature leads to faster disintegration compared to random removal of edges in the model and real directed networks (Fig. \ref{rob_edge}). Moreover, removal of edges based on increasing order of FR or AFR curvature leads to slightly faster disintegration in comparison to removal of edges based on increasing order of OR curvature, in model and real directed networks considered here (Fig. \ref{rob_edge}). However, removal of edges based on decreasing order of edge betweenness centrality typically leads to slightly faster disintegration in comparison to the removal of edges based on increasing order of OR, FR or AFR curvature in model and real directed networks considered here (Fig. \ref{rob_edge}). These results underscore that our definition of the OR, FR and AFR curvature captures the importance of edges for flows in directed networks.

\begin{figure}
\centering
\includegraphics[width=.97\columnwidth]{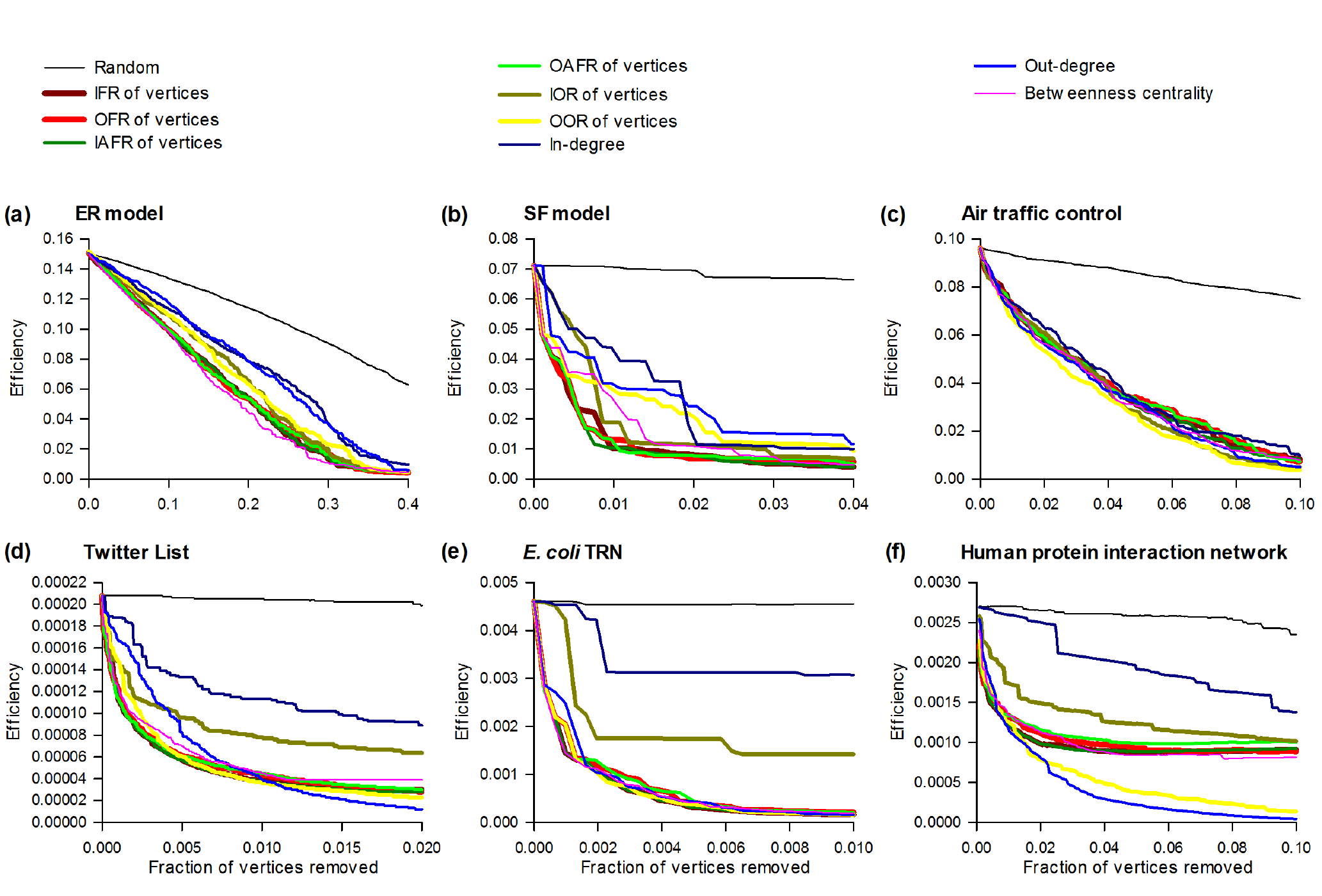}
\caption{Communication efficiency as a function of the fraction of vertices removed in model and real-world directed networks. (a) ER model with $n=1000$, $p=0.003$. (b) SF model with $n=1000$, $m=3000$, $\lambda_{in}$=2.1, $\lambda_{out}$=2.1. (c) Air traffic control. (d) Twitter List. (e)  \textit{E. coli} TRN. (f) Human protein interaction network.}
\label{rob_vertex}
\end{figure}

We next investigated the effect of removing vertices based on increasing order of their IOR, OOR, IFR, OFR, IAFR or OAFR curvature on the large-scale connectivity of directed networks (Fig. \ref{rob_vertex}). It is seen that targeted removal of vertices with highly negative IOR, OOR, IFR, OFR, IAFR or OAFR curvature leads to faster disintegration compared to random removal of vertices in both model and real directed networks (Fig. \ref{rob_vertex}). Previous work \cite{Barabasi1999,Jeong2000,Jeong2001,Albert2002,Joy2005,Yu2007,Newman2010} has shown that model and real networks are vulnerable to targeted removal of vertices with high degree or high betweenness centrality. We have also compared the effect of removing vertices based on increasing order of their IOR, OOR, IFR, OFR, IAFR or OAFR curvature against removing vertices based on decreasing order of their in-degree, out-degree or betweenness centrality on the large-scale connectivity of networks (Fig. \ref{rob_vertex}). In majority of model and real-world directed networks considered here, we find that the removal of vertices based on decreasing order of betweenness centrality leads to at least slightly faster distintegration compared to removal of vertices based on increasing order of IOR, OOR, IFR, OFR, IAFR or OAFR curvature or decreasing order of in-degree or out-degree (Fig. \ref{rob_vertex}). These results underscore that vertices with highly negative IOR, OOR, IFR, OFR, IAFR or OAFR curvature are important for maintaining the large-scale connectivity of directed networks.



\section{Summary and Future Outlook}
\label{summary}

In this paper, we have extended both our recent adaptation of the Forman's discretization of Ricci curvature and its extension to include two-dimensional faces capturing higher-order correlations between vertices, as well as Ollivier's discretization of Ricci curvature, from undirected to directed networks.

The mathematical expression for the Forman-Ricci (FR) curvature for a directed edge elegantly incorporates the vertex weights, edge weights and edge direction in directed networks. Furthermore, the augmented Forman-Ricci (AFR) curvature takes into account also the influence of the faces adjacent to a given directed edge on the Ricci curvature of the considered edge. We find that the distribution of the FR and AFR curvature for directed edges is narrow in directed ER networks, while the distribution is broad in directed SF networks. In most real directed networks, the distribution of the FR and AFR curvature for directed edges is also broad like SF networks. These results highlight that the distribution of the FR and AFR curvature for directed edges can also be employed to distinguish and classify different types of directed networks. We also find that the FR and AFR curvature for edges has a moderate to high negative correlation with edge betweenness centrality in considered directed networks. Based on the definition of the FR and AFR curvature for a directed edge, it is easy to define the In Forman-Ricci (IFR), Out Forman-Ricci (OFR), In Augmented Forman-Ricci (IAFR) and Out Augmented Forman-Ricci (OAFR) curvature for vertices in directed networks. We have also investigated the correlation of IFR and IAFR curvature with in-degree, OFR and OAFR curvature with out-degree, and IFR, OFR, IAFR and OAFR curvature with betweenness centrality and page rank of vertices in model and real directed networks. Moreover, an investigation of the effect of removing edges based on their FR or AFR curvature, or removing vertices based on their IFR, OFR, IAFR or OAFR curvature, on the communication efficiency of directed networks shows that both model and real directed networks are vulnerable to targeted removal of edges with highly negative FR or AFR curvature or vertices with highly negative IFR, OFR, IAFR or OAFR curvature.

Using the same basic approach employed to extend Forman's discretization of Ricci curvature for undirected networks to directed networks, we have extended the better known Ollivier-Ricci (OR) curvature for undirected networks to directed networks (See section \ref{theory}). We find that the OR curvature for directed edges has a high negative correlation with FR or AFR curvature in many model and real directed networks analyzed here. We emphasize that, in contrast to the FR or AFR curvature whose computation is simple and fast in large undirected and directed networks, the computation of the OR curvature necessitates solving a linear programming problem associated with optimal mass transport on networks which will not computationally scale to extremely large directed networks such as the world wide web (WWW) \cite{Broder2000}. Hence, our result that OR curvature for directed edges has a high negative correlation with FR or AFR curvature in many real directed networks suggests that FR or AFR curvature can be used to obtain a coarse-view on OR curvature in extremely large directed networks where the computation of OR curvature may be prohibitive. We also find that the OR curvature for edges has a moderate to high negative correlation with edge betweenness centrality in considered directed networks. We have also investigated the correlation of IOR curvature with in-degree, OOR curvature with out-degree, and IOR and OOR curvature with betweenness centrality and page rank of vertices in model and real directed networks. Also, an investigation of the effect of removing edges based on their OR curvature, or removing vertices based on their IOR or OOR curvature, on the communication efficiency of directed networks shows that both model and real directed networks are vulnerable to targeted removal of edges with highly negative OR curvature or vertices with highly negative IOR or OOR curvature. Thus, both Forman's and Ollivier's discretizations of Ricci curvature can hereafter be employed to uncover the hidden geometry of both directed and undirected complex networks.

Notably, a drawback of the FR and OR curvature is that they do not take into account higher order correlations between triples, quadruples and higher n-tuples of vertices. However, we could overcome this drawback by appealing to Forman's original definition of Ricci curvature \cite{Forman2003}, which is defined for polyhedral and more general complexes, and thus, intrinsically takes into account faces of dimension higher or equal to two which can capture the higher-order correlations between vertices in networks \cite{Sreejith2016,Weber2017,Samal2018}. Importantly, we have also introduced here the adaptation of the augmented Forman-Ricci (AFR) curvature for undirected networks to directed networks which account for the two-dimensional directed simplicial complexes. Among the three essentially different types of directed triangles (Fig. \ref{DirectedExplain}(b)), we have chosen the feed forward loop (FFL) to be the face or directed simplicial complex while computing the AFR curvature for directed edges (Fig. \ref{DirectedExplain}(c). Note that more directed edges in analyzed model and real directed networks have a positive AFR curvature in comparison to positive FR curvature, and this is presumably due to positive contributions from directed simplicial complexes accounted in the definition of the AFR curvature.

A number of potential subjects of further research present themselves as natural extensions of this work. Perhaps the first and most basic among them is the exploration of the augmented Ricci flow in different types of networks. Previously, some of us have explored the Ricci flow in undirected networks \cite{Weber2017}, and have further shown the applications of Ricci flow in the prediction of the long time behavior of the network. Clearly, a similar and systematic investigation of the Ricci flow with the definition of the AFR curvature for directed edges is warranted in model and real-world directed networks. Moreover, the long time evolution is given by the Euler-Bloch characteristic of the network, and a similar invariant also needs to be developed for directed networks. These research directions represent work in progress and will be reported in future contributions. In a nutshell, we expect the definitions of the discrete Ricci curvatures introduced here, especially the AFR curvature, to be employed for the future investigations of the directed networks and hyper-networks.


\section*{Acknowledgments}
We thank M. Karthikeyan for help with figures. E.S. and A.S. thank the Max Planck Institute for Mathematics in the Sciences, Leipzig, for their warm hospitality. A.S. would like to acknowledge financial support from Max Planck Society through the award of a Max Planck Partner Group in Mathematical Biology.



\begin{table}
\caption{\label{tab-1}Comparison of the Ollivier-Ricci (OR) with Forman-Ricci (FR) and Augmented Forman-Ricci (AFR) curvature for edges in model and real directed networks. In this table, we report the Spearman correlation of OR with FR or AFR curvature for edges. In case of model networks, the reported correlation is mean (rounded to two decimal places) and standard error (rounded to four decimal places) over a sample of 100 directed networks generated with specific input parameters.}
\centering
\begin{tabular}{|l|c|c|}
\hline
\multicolumn{1}{|c|}{\small Network} & \multicolumn{2}{|c|}{\small OR versus} \\
\cline{2-3}
     & {\small FR} & {\small AFR} \\
\hline
\multicolumn{3}{|c|}{\small Model networks} \\
\hline
{\small ER model with n=1000, p=0.003} & {\small 0.82 $\pm$ 0.0117} & {\small 0.82 $\pm$ 0.0117} \\
\hline
{\small ER model with n=1000, p=0.005} & {\small 0.44 $\pm$ 0.0238} & {\small 0.44 $\pm$ 0.0238} \\
\hline
{\small ER model with n=1000, p=0.01} & {\small -0.28 $\pm$ 0.0148} & {\small -0.28 $\pm$ 0.0148} \\
\hline
{\small ER model with n=1000, p=0.02} & {\small -0.37 $\pm$ 0.0096} & {\small -0.37 $\pm$ 0.0096} \\
\hline
{\small SF model with n=1000, m=3000, $\lambda_{in}$=2.1, $\lambda_{out}$=2.1} & {\small 0.59 $\pm$ 0.0355} & {\small 0.58 $\pm$ 0.0337} \\
\hline
{\small SF model with n=1000, m=5000, $\lambda_{in}$=2.1, $\lambda_{out}$=2.1} & {\small 0.32 $\pm$ 0.0469} & {\small 0.33 $\pm$ 0.0477} \\
\hline
{\small SF model with n=1000, m=10000, $\lambda_{in}$=2.1, $\lambda_{out}$=2.1} & {\small 0.02 $\pm$ 0.0441} & {\small 0.05 $\pm$ 0.0486} \\
\hline
{\small SF model with n=1000, m=20000, $\lambda_{in}$=2.1, $\lambda_{out}$=2.1} & {\small -0.12 $\pm$ 0.0393} & {\small -0.08 $\pm$ 0.0437} \\
\hline
\multicolumn{3}{|c|}{\small Real networks} \\
\hline
{\small Air traffic control} & {\small 0.29}  & {\small 0.30} \\
\hline
{\small Twitter List} & {\small 0.87} & {\small 0.86} \\
\hline
{\small \textit{E. coli} TRN} & {\small 0.53} & {\small 0.42} \\
\hline
{\small \textit{B. subtilis} TRN} & {\small 0.83} & {\small 0.74} \\
\hline
{\small Human protein interaction network} & {\small 0.97} & {\small 0.92} \\
\hline
{\small Phosphonetwork} & {\small 0.95} & {\small 0.90} \\
\hline
\end{tabular}
\end{table}

\pagestyle{empty}
\begin{landscape}
\begin{table}
\caption{\label{tab-2}Comparison of edge betweenness centrality (EBC) with Ollivier-Ricci (OR), Forman-Ricci (FR) and Augmented Forman-Ricci (AFR) curvature for edges in model and real directed networks. In this table, we report the Spearman correlation of EBC with OR, FR and AFR curvature for edges. In case of model networks, the reported correlation is mean (rounded to two decimal places) and standard error (rounded to four decimal places) over a sample of 100 directed networks generated with specific input parameters.}
\centering
\begin{tabular}{|l|c|c|c|}
\hline
\multicolumn{1}{|c|}{\small Network} & \multicolumn{3}{|c|}{\small EBC versus} \\
\cline{2-4}
         & {\small OR} & {\small FR} & {\small AFR} \\
\hline
\multicolumn{4}{|c|}{\small Model networks} \\
\hline
{\small ER model with n=1000, p=0.003} & {\small -0.74 $\pm$ 0.0138} & {\small  -0.77 $\pm$ 0.0119} & {\small -0.77 $\pm$ 0.0118} \\
\hline
{\small ER model with n=1000, p=0.005} & {\small -0.46 $\pm$ 0.0222} & {\small  -0.81 $\pm$ 0.0084} & {\small  -0.81 $\pm$ 0.0083} \\
\hline
{\small ER model with n=1000, p=0.01} & {\small 0.17 $\pm$ 0.0137} & {\small  -0.83 $\pm$ 0.0071} & {\small -0.83 $\pm$ 0.0068} \\
\hline
{\small ER model with n=1000, p=0.02} & {\small 0.26 $\pm$ 0.0105} & {\small  -0.79 $\pm$ 0.0081} & {\small  -0.80 $\pm$ 0.0079} \\
\hline
{\small SF model with n=1000, m=3000, $\lambda_{in}$=2.1, $\lambda_{out}$=2.1} & {\small -0.63 $\pm$ 0.0568} & {\small  -0.66 $\pm$ 0.0386} & {\small  -0.69 $\pm$ 0.0360} \\
\hline
{\small SF model with n=1000, m=5000, $\lambda_{in}$=2.1, $\lambda_{out}$=2.1} & {\small -0.42 $\pm$ 0.0707} & {\small  -0.57 $\pm$ 0.0354} & {\small  -0.64 $\pm$ 0.0295} \\
\hline
{\small SF model with n=1000, m=10000, $\lambda_{in}$=2.1, $\lambda_{out}$=2.1} & {\small -0.18 $\pm$ 0.0591} & {\small  -0.57 $\pm$ 0.0291} & {\small -0.68 $\pm$ 0.0269} \\
\hline
{\small SF model with n=1000, m=20000, $\lambda_{in}$=2.1, $\lambda_{out}$=2.1} & {\small 0.16 $\pm$ 0.4218} & {\small -0.60 $\pm$ 0.0257} & {\small -0.73 $\pm$ 0.0173} \\
\hline
\multicolumn{4}{|c|}{\small Real networks} \\
\hline
{\small Air traffic control} & {\small -0.50} & {\small -0.38} & {\small -0.39} \\
\hline
{\small Twitter List} & {\small -0.13} & {\small -0.22} & {\small -0.22} \\
\hline
{\small \textit{E. coli} TRN} & {\small -0.70} & {\small -0.49} & {\small -0.50} \\
\hline
{\small \textit{B. subtilis} TRN} & {\small -0.83} & {\small -0.73} & {\small -0.70} \\
\hline
{\small Human protein interaction network} & {\small -0.83} & {\small -0.85} & {\small -0.82} \\
\hline
{\small Phosphonetwork} & {\small -0.84} & {\small -0.83} & {\small -0.82} \\
\hline
\end{tabular}
\end{table}
\end{landscape}

\pagestyle{empty}
\begin{landscape}
\begin{table}
\caption{\label{tab-3}Comparison of In Ollivier-Ricci (IOR) with In Forman-Ricci (IFR) and In Augmented Forman-Ricci (IAFR) curvature for vertices, and Out Olliver-Ricci (OOR) with Out Forman-Ricci (OFR) and Out Augmented Forman-Ricci (OAFR) curvature for vertices in model and real directed networks. In this table, we report the Spearman correlation of IOR with IFR and IAFR curvature for vertices, and OOR with OFR and OAFR curvature for vertices. In case of model networks, the reported correlation is mean (rounded to two decimal places) and standard error (rounded to four decimal places) over a sample of 100 directed networks generated with specific input parameters.}
\centering
\begin{tabular}{|l|c|c|c|c|}
\hline
\multicolumn{1}{|c|}{\small Network} & \multicolumn{2}{|c|}{\small IOR versus} & \multicolumn{2}{|c|}{\small OOR versus} \\
\cline{2-5}
  & {\small IFR} & {\small IAFR} & {\small OFR} & {\small OAFR} \\
\hline
\multicolumn{5}{|c|}{\small Model networks} \\
\hline
{\small ER model with n=1000, p=0.003} & {\small 0.87 $\pm$ 0.0111} & {\small 0.87 $\pm$ 0.0112} & {\small 0.87 $\pm$ 0.0087} & {\small 0.87 $\pm$ 0.0088} \\
\hline
{\small ER model with n=1000, p=0.005} & {\small 0.86 $\pm$ 0.0097} & {\small 0.86 $\pm$ 0.0097} & {\small 0.86 $\pm$ 0.0089} & {\small 0.86 $\pm$ 0.0090} \\
\hline
{\small ER model with n=1000, p=0.01} & {\small 0.79 $\pm$ 0.0117} & {\small 0.79 $\pm$ 0.0119} & {\small 0.79 $\pm$ 0.0136} & {\small 0.79 $\pm$ 0.0138} \\
\hline
{\small ER model with n=1000, p=0.02} & {\small 0.75 $\pm$ 0.0154} & {\small 0.75 $\pm$ 0.0157} & {\small 0.75 $\pm$ 0.0162} & {\small 0.75 $\pm$ 0.0165} \\
\hline
{\small SF model with n=1000, m=3000, $\lambda_{in}$=2.1, $\lambda_{out}$=2.1} & {\small 0.78 $\pm$ 0.0223} & {\small 0.76 $\pm$ 0.0252} & {\small 0.78 $\pm$ 0.0247} & {\small 0.76 $\pm$ 0.0282} \\
\hline
{\small SF model with n=1000, m=5000, $\lambda_{in}$=2.1, $\lambda_{out}$=2.1} & {\small 0.77 $\pm$ 0.0187} & {\small 0.75 $\pm$ 0.0210} & {\small 0.77 $\pm$ 0.0181} & {\small 0.75 $\pm$ 0.0212} \\
\hline
{\small SF model with n=1000, m=10000, $\lambda_{in}$=2.1, $\lambda_{out}$=2.1} & {\small 0.72 $\pm$ 0.0171} & {\small 0.71 $\pm$ 0.0180} & {\small 0.72 $\pm$ 0.0168} & {\small 0.71 $\pm$ 0.0174} \\
\hline
{\small SF model with n=1000, m=20000, $\lambda_{in}$=2.1, $\lambda_{out}$=2.1} & {\small 0.66 $\pm$ 0.0157} & {\small 0.64 $\pm$ 0.0160} & {\small 0.66 $\pm$ 0.0181} & {\small 0.64 $\pm$ 0.0186} \\
\hline
\multicolumn{5}{|c|}{\small Real networks} \\
\hline
{\small Air traffic control} & {\small 0.76} & {\small 0.76} & {\small 0.86} & {\small 0.86} \\
\hline
{\small Twitter List} & {\small 0.84} & {\small 0.77} & {\small 0.51} & {\small 0.50} \\
\hline
{\small \textit{E. coli} TRN} & {\small 0.74} & {\small 0.30} & {\small 0.33} & {\small 0.33} \\
\hline
{\small \textit{B. subtilis} TRN} & {\small 0.65} & {\small 0.50} & {\small -0.08} & {\small -0.08} \\
\hline
{\small Human protein interaction network} & {\small 0.29} & {\small 0.23} & {\small 0.48} & {\small 0.47} \\
\hline
{\small Phosphonetwork} & {\small 0.72} & {\small 0.69} & {\small 0.55} & {\small 0.55} \\
\hline
\end{tabular}
\end{table}
\end{landscape}

\pagestyle{empty}
\begin{landscape}
\begin{table}
\caption{\label{tab-4}Comparison of in-degree with In Ollivier-Ricci (IOR), In Forman-Ricci (IFR) and In Augmented Forman-Ricci (IAFR) curvature for vertices, and out-degree with Out Olliver-Ricci (OOR), Out Forman-Ricci (OFR) and Out Augmented Forman-Ricci (OAFR) curvature for vertices in model and real directed networks. In this table, we report the Spearman correlation of in-degree with IOR, IFR and IAFR curvature for vertices, and out-degree with OOR, OFR and OAFR curvature for vertices. In case of model networks, the reported correlation is mean (rounded to two decimal places) and standard error (rounded to four decimal places) over a sample of 100 directed networks generated with specific input parameters.}
\centering
\begin{tabular}{|l|c|c|c|c|c|c|}
\hline
\multicolumn{1}{|c|}{\small Network} & \multicolumn{3}{|c|}{\small In-degree versus} & \multicolumn{3}{|c|}{\small Out-degree versus} \\
\cline{2-7}
  & {\small IOR} & {\small IFR} & {\small IAFR} & {\small OOR} & {\small OFR} & {\small OAFR} \\
\hline
\multicolumn{7}{|c|}{\small Model networks} \\
\hline
{\small ER model with n=1000, p=0.003} & {\small -0.95 $\pm$ 0.0052} & {\small -0.73 $\pm$ 0.0199} & {\small -0.73 $\pm$ 0.0201} & {\small -0.95 $\pm$ 0.0060} & {\small -0.74 $\pm$ 0.0175} & {\small -0.73 $\pm$ 0.0176} \\
\hline
{\small ER model with n=1000, p=0.005} & {\small -0.98 $\pm$ 0.0029} & {\small -0.81 $\pm$ 0.0125} & {\small -0.81 $\pm$ 0.0126} & {\small -0.98 $\pm$ 0.0028} & {\small -0.81 $\pm$ 0.0133} & {\small -0.81 $\pm$ 0.0134} \\
\hline
{\small ER model with n=1000, p=0.01} & {\small -0.98 $\pm$ 0.0014} & {\small -0.85 $\pm$ 0.0093} & {\small -0.85 $\pm$ 0.0095} & {\small -0.98 $\pm$ 0.0014} & {\small -0.85 $\pm$ 0.0102} & {\small -0.85 $\pm$ 0.0104} \\
\hline
{\small ER model with n=1000, p=0.02} & {\small -0.97 $\pm$ 0.0026} & {\small -0.87 $\pm$ 0.0084} & {\small -0.87 $\pm$ 0.0087} & {\small -0.97 $\pm$ 0.0028} & {\small -0.87 $\pm$ 0.0092} & {\small -0.87 $\pm$ 0.0095} \\
\hline
{\small SF model with n=1000, m=3000, $\lambda_{in}$=2.1, $\lambda_{out}$=2.1} & {\small -0.88 $\pm$ 0.0188} & {\small -0.70 $\pm$ 0.0332} & {\small -0.67 $\pm$ 0.0357} & {\small -0.88 $\pm$ 0.0192} & {\small -0.70 $\pm$ 0.0343} & {\small -0.67 $\pm$ 0.0374} \\
\hline
{\small SF model with n=1000, m=5000, $\lambda_{in}$=2.1, $\lambda_{out}$=2.1} & {\small -0.92 $\pm$ 0.0094} & {\small -0.75 $\pm$ 0.0233} & {\small -0.73 $\pm$ 0.0264} & {\small -0.92 $\pm$ 0.0097} & {\small -0.75 $\pm$ 0.0232} & {\small -0.73 $\pm$ 0.0270} \\
\hline
{\small SF model with n=1000, m=10000, $\lambda_{in}$=2.1, $\lambda_{out}$=2.1} & {\small -0.90 $\pm$ 0.0103} & {\small -0.79 $\pm$ 0.0184} & {\small -0.77 $\pm$ 0.0197} & {\small -0.91 $\pm$ 0.0074} & {\small -0.79 $\pm$ 0.0180} & {\small -0.77 $\pm$ 0.0188} \\
\hline
{\small SF model with n=1000, m=20000, $\lambda_{in}$=2.1, $\lambda_{out}$=2.1} & {\small -0.85 $\pm$ 0.0125} & {\small -0.81 $\pm$ 0.0181} & {\small -0.79 $\pm$ 0.0188} & {\small -0.85 $\pm$ 0.0117} & {\small -0.81 $\pm$ 0.0169} & {\small -0.79 $\pm$ 0.0173} \\
\hline
\multicolumn{7}{|c|}{\small Real networks} \\
\hline
{\small Air traffic control} & {\small -0.84} & {\small -0.69} & {\small -0.68} & {\small -0.89} & {\small -0.78} & {\small -0.78} \\
\hline
{\small Twitter List} & {\small -0.35} & {\small -0.17} & {\small -0.01} & {\small -0.51} & {\small -0.90} & {\small -0.89} \\
\hline
{\small \textit{E. coli} TRN} & {\small -0.83} & {\small -0.46} & {\small 0.10} & {\small -0.33} & {\small -0.99} & {\small -0.99} \\
\hline
{\small \textit{B. subtilis} TRN} & {\small -0.47} & {\small -0.02} & {\small 0.29} & {\small 0.07} & {\small -0.95} & {\small -0.95} \\
\hline
{\small Human protein interaction network} & {\small -0.49} & {\small 0.08} & {\small 0.13} & {\small -0.50} & {\small -0.94} & {\small -0.93} \\
\hline
{\small Phosphonetwork} & {\small -0.70} & {\small -0.44} & {\small -0.41} & {\small -0.55} & {\small -0.98} & {\small -0.98} \\
\hline
\end{tabular}
\end{table}
\end{landscape}

\newpage
\section*{Supplementary Information}

\noindent \textbf{Supplementary Table S1:} Description of considered model and real-world directed networks. The table lists the number of vertices, number of edges, maximum degree, average degree, average In-degree, average Out-degree and size of the largest weakly connected component for each directed network.

\noindent \textbf{Supplementary Table S2:} Comparison of Betweenness centrality with In Ollivier-Ricc (IOR), Out Ollivier-Ricci (OOR), In Forman-Ricci (IFR), Out Forman-Ricci (OFR), In Augmented Forman-Ricci (IAFR) and Out Augmented Forman-Ricci (OAFR) curvature for vertices in model and real directed networks. In this table, we report the Spearman correlation of Betweenness centrality with IOR, OOR, IFR, OFR, IAFR and OAFR for vertices. In case of model networks, the reported correlation is mean (rounded to two decimal places) and standard error (rounded to four decimal places) over a sample of 100 directed networks generated with specific input parameters.

\noindent \textbf{Supplementary Table S3:} Comparison of Page rank with In Ollivier-Ricci (IOR), Out Ollivier-Ricci (OOR), In Forman-Ricci (IFR), Out Forman-Ricci (OFR), In Augmented Forman-Ricci (IAFR) and Out Augmented Forman-Ricci (OAFR) curvature for vertices in model and real directed networks. In this table, we report the Spearman correlation of Page rank with IOR, OOR, IFR, OFR, IAFR and OAFR for vertices. In case of model networks, the reported correlation is mean (rounded to two decimal places) and standard error (rounded to four decimal places) over a sample of 100 directed networks generated with specific input parameters.	


\end{document}